\documentclass[article]{amsart}
\usepackage[textheight=8.5in,margin=1.6in]{geometry}
\usepackage[foot]{amsaddr}
\usepackage{amsmath}
\usepackage{amsfonts}
\usepackage{graphicx}
\usepackage{epstopdf}
\usepackage{url}

\usepackage{algorithmic,algorithm}
\ifpdf
\DeclareGraphicsExtensions{.eps,.pdf,.png,.jpg}
\else
\DeclareGraphicsExtensions{.eps}
\fi

\usepackage{amsopn}

\newcommand{\norm}[1]{\left\lVert#1\right\rVert}
\newcommand{\abs}[1]{\left\lvert#1\right\rvert}

\DeclareMathOperator*{\argmin}{arg\,min}

\newcommand{\R}{\mathbb{R}}
\newcommand{\N}{\mathbb{N}}

\newtheorem*{remark}{Remark}
\newtheorem{theorem}{Theorem}
\newtheorem{lemma}{Lemma}

\begin{document}

\title[Fast Flexible LSQR for Large-Scale Regularization]{Fast Flexible LSQR with a Hybrid Variant for Efficient Large-Scale Regularization}	
\thanks{The work of the authors was partially supported by the Charles University, project GAUK 376121, and partially by the grant SVV-2023-260711.}

\author[Eva Miku\v{s}ov{\'a}]{ Eva Miku\v{s}ov{\'a}}
\email{eva.havelkova@karlin.mff.cuni.cz (Eva Mikušová)}

\author[Iveta Hn{\v{e}}tynkov{\'a}]{Iveta Hn{\v{e}}tynkov{\'a}}
\email{iveta.hnetynkova@karlin.mff.cuni.cz (Iveta Hnětynková)}

\address{Department of Numerical Mathematics, Faculty of Mathematics and Physics, Charles University, Sokolovská 83, 18675 Prague, Czech Republic}



\begin{abstract}
A wide range of applications necessitates solving large-scale ill-posed problems
 contaminated by noise. Krylov subspace regularization methods are particularly advantageous 
in this context, as they rely solely on matrix-vector multiplication.  
Among the most widely used techniques are LSQR and CGLS, both of which can be extended with 
flexible preconditioning to enforce solution properties such as nonnegativity or sparsity. 
Flexible LSQR (FLSQR) can also be further combined with direct methods to create efficient hybrid approaches. 
The Flexible Golub-Kahan bidiagonalization underlying FLSQR requires two long-term recurrences. 
In this paper, we introduce a novel Fast Flexible Golub-Kahan bidiagonalization method that 
employs one long-term and one short-term recurrence. Using this, we develop the Fast Flexible 
LSQR (FaFLSQR) algorithm, which offers comparable computational cost to FCGLS while 
also supporting hybrid regularization like FLSQR. We analyze the properties of FaFLSQR and prove 
its mathematical equivalence to FCGLS. Numerical experiments demonstrate that in floating point arithmetic
FaFLSQR outperforms both FCGLS and FLSQR in terms of computational efficiency. 
 
\end{abstract}

\maketitle
\section{Introduction}
We are interested in solving linear ill-posed inverse problems of the form
\begin{equation}\label{theProblem}
	A x \approx b, \quad A\in \R^{m \times n}, b \in \R^m,
\end{equation}
where the vector $b$ represents the observed noisy data and $A$ is an ill-conditioned matrix 
representing the model. In real-world applications, the matrix $A$ is typically large scale 
and may not be given explicitly. Furthermore, its singular values decay gradually to zero without a
clear gap, so its numerically rank is not well defined. Consequently, multiplication by $A$ smooths
out high frequencies in a vector (the so called smoothing property).  

The observation $b$ contains noise (errors) from various sources. We assume that $b=b^{exact}+e$, where $e$ 
is some unknown additive white noise, and $b^{exact}$ is the unknown exact observation. Our goal is to approximate
$x^{exact}$ such that $A x^{exact} = b^{exact}$, while only having access to the data in the equation \eqref{theProblem}.
Denote $\eta=\norm{e}/\norm{b^{exact}}$ the noise level in the data.

It is well know that attempt to solve \eqref{theProblem} directly fails, as noise tends to get amplified 
during the process. Therefore, regularization is necessary. When $A$ is large or not explicitly given, 
iterative Krylov subspace methods that rely on matrix-vector multiplication offer a distinct advantage.
In these methods, the original problem is iteratively projected onto a sequence of Krylov subspaces with increasing dimension, and the solution is approximated based on a small, reduced least squares subproblem. The regularization effect arises from the fact that Krylov subspaces are in early iterations dominated by smooth vectors, which tend to suppress high-frequency noise. However, as the number of iterations increases, noise begins to propagate into the projection, a phenomenon known as semiconvergence. This effect necessitates early stopping of the iterations to prevent overfitting to the noisy data. In this context, the number of iterations effectively serves as the regularization parameter.

Two powerful methods for solving large inverse problems with rectangular model matrices are LSQR (Least Squares QR)~\cite{LSQR, LSQR2} and CGLS (Conjugate Gradient Least Squares)~\cite{CGLS}. While mathematically equivalent, these methods work differently and each has its own advantages. Both can be extended with flexible (iteration-dependent) preconditioning to enforce specific solution properties, such as sparsity or nonnegativity, further enhancing their practical utility in a variety of applications, see, e.g., \cite{flexibleKrylovLp, iterativelyReweighted, nonNegative, flexibleKrylovAugment}. 
LSQR, in its basic form, relies on the Golub-Kahan bidiagonalization (GK) process, which constructs two sets of orthonormal vectors using short recurrences. To incorporate a flexible preconditioner, the GK process was replaced in \cite{flexibleKrylovLp} by the Flexible Golub-Kahan (FGK) process, which requires two long recurrences. This increases the computational cost of the resulting Flexible LSQR (FLSQR). Flexible CGLS (FCGLS) \cite{nonNegative} requires only one long recurrence but does not allow for a hybrid framework. In the hybrid version of Krylov subspace methods, the projected subproblem is in each outer iteration solved using a direct regularization method, such as Tikhonov regularization; see, e.g., \cite{flexibleKrylovLp, EvickaIvetka, RegParComparison, RenautHybrid, RenautHybridRew}
and also \cite{surveyHybrid} for an overview. If the inner regularization parameter is chosen appropriately, this approach helps to overcome the semiconvergence issue. The solution is stabilized after a certain number of iterations, allowing for an easy stopping. This effectively means that the regularization parameter selection is shifted to the smaller projected problem.

{\em Key contribution:}
In this paper, we introduce the Fast Flexible LSQR (FaFLSQR) method based on a novel version of the Golub-Kahan iterative bidiagonalization. Our Fast Flexible Golub-Kahan process requires one long and one short recurrence to generate two sequences of basis vectors with particular orthogonality properties. This makes FaFLSQR computationally more efficient than the existing FLSQR. We prove that, in exact arithmetic, FaFLSQR is mathematically equivalent to FCGLS and it has comparable computational cost per iteration. Furthermore, FaFLSQR allows for the incorporation of inner regularization within a hybrid framework. Thus, FaLSQR has computational efficiency of FCGLS and allows extensions as FLSQR.
Finally, we demonstrate through numerical experiments that the novel FaLSQR outperforms not only FLSQR but also FCGLS in terms of computational efficiency. This is because in floating point arithmetic FaFLSQR does not suffer from the delay of convergence as much as FCGLS.

The paper is organized as follows. Section 2 summarizes CGLS and LSQR algorithms and their known flexible variants. Our main results are given in section 3 presenting the Fast Flexible Golub-Kahan bidiagonalization and Fast Flexible LSQR. Hybrid version is discussed. Section 4 compares the new algorithm to FLSQR, and proofs its equivalency to FCGLS. Computational cost is discussed. Section 5 contains numerical experiments showing the performance of the new algorithm in comparison to other flexible methods. The conclusions follow in section 6. Auxiliary proofs are moved to the appendix. For simplicity of exposition we assume in the whole paper that the algorithms do not stop before the considered iterations.

\section{Preconditioning of CGLS and LSQR}
CGLS and LSQR are mathematically equi\-va\-lent iterative methods developed originally for the solution 
of the least squares problem within \eqref{theProblem}. Consider the same initial approximation $x_0$ for all algorithms in this paper.
Let $r_0 = b - A x_0$ be the corresponding initial residual. CGLS and LSQR are both based on short recurrences and search the $k$-th approximate solution in the form $\bar{x}_k \in x_0 + \mathcal{K}_k(A^TA, A^Tr_0)$,
where $\mathcal{K}_k(A^TA, A^Tr_0)$ is the Krylov subspace. The approximation fulfills the optimality conditions
\begin{equation}\label{lsqrOriginalProblem}
\bar{x}_k=\argmin_{x \in x_0 + \mathcal{K}_k(A^TA, A^Tr_0)}{\norm{x_{LS}-x}_{A^TA}} 
= \argmin_{x \in x_0 + \mathcal{K}_k(A^TA, A^Tr_0)}{\norm{b-Ax}},
\end{equation}
where $x_{LS}$ is a least squares solution of \eqref{theProblem}.

CGLS was derived by applying CG method \cite{CGLS} to normal equations $A^TA x \approx A^Tb$ and rearranged to avoid work
with the matrix $A^TA$. It generates in $k$ iterations a set of direction vectors $\bar{p}_0, \dots, \bar{p}_{k-1} \in \R^{n}$ representing
an $A^TA$-orthogonal basis of $\mathcal{K}_k (A^TA, A^Tr_0)$. 
The approximation $\bar{x}_k$ is then updated iteratively in the selected directions.
LSQR was derived independently based on the Golub-Kahan iterative bidiagonalization process \cite{LSQR, LSQR2}. In iteration $k$ 
it constructs matrices $\bar{V}_k = [\bar{v}_1, \dots, \bar{v}_k] \in \R^{n \times k}$ and $\bar{U}_k= [\bar{u}_1, \dots, \bar{u}_k ] \in \R^{m \times k}$ 
with columns representing orthonormal bases of $\mathcal{K}_k(A^TA,A^Tr_0)$
and $\mathcal{K}_k(AA^T,r_0)$, respectively. The bidiagonalization (orthogonalization) coefficients are saved
in the lower bidiagonal matrix $\bar{L}_k \in \R^{k \times k}$ and the following relations holds
\[
	A\bar{V}_k =\bar{U}_{k+1}\bar{L}_{k+}, \qquad A^T\bar{U}_k =\bar{V}_k\bar{L}_k^T,
\]
where $\bar{L}_{k+} \in \R^{(k+1) \times k}$ is the matrix $\bar{L}_k$ extended by one more row. The approximation $\bar{x}_k$ has then the form $\bar{x}_k = x_0+\bar{V}_k \bar{y}_k$, where $\bar{y}_k$ is the least squares solution of the small projected problem $\bar{L}_{k+} \bar{y}_k \approx \beta_1 e_1$, $\beta_1 =||r_0||$. 

A variety of preconditioning strategies are available for both LSQR and CGLS. One common approach is to use a suitable preconditioning matrix  $B$ such that $BB^T \approx A^T A.$ For LSQR, right preconditioning is typically applied, and the analytically equivalent strategy for CGLS is known as split preconditioning; see, e.g., \cite{Greenbaum, HansenRankDef, Bjorck} for more references. However, in some cases an explicit factorization $BB^T \approx A^T A$ may not be available or may be too costly to compute. Instead, only a preconditioner $M \approx A^T A $ may be accessible. For such situations, factorization-free variants of both CGLS and LSQR have been developed, for further information see \cite{PCGLS, saad, ArridgePreconditionedLSQR}. We provide a more detailed discussion of the latter at the beginning of section~\ref{Section3}.  For additional analysis, properties, and theoretical connections, we refer the reader to~\cite{EvickaIvetkaEnumath}.

In flexible (iteration dependent) preconditioning, the preconditioner is allowed to vary at each iteration. Consider a sequence of symmetric positive-definite preconditioners $M_k, k=1, 2, ...$ (often diagonal in practice). Incorporating flexible preconditioning results in a significant change in the properties of algorithms. Specifically, short recurrences fail to produce orthogonal bases for the desired Krylov subspaces, necessitating the use of full recurrences. In the Flexible Golub-Kahan bidiagonalization (FGK) presented in 
\cite{flexibleKrylovLp}, both the 3-term recurrences used to 
generate the vectors $\bar{u}_k$ and $\bar{v}_k$ are replaced with full recurrences, as seen in lines 10 and 17 of 
Algorithm~\ref{AlgorithmFGK}. Then $\tilde{V}_k = [\tilde{v}_1, \dots, \tilde{v}_k] $ and $\tilde{U}_k = [\tilde{u}_1, \dots, \tilde{u}_k]$ have orthonormal columns. For the outputs, it holds that
\begin{equation}\label{FGKdecomposition}	
		A\tilde{Z}_k=\tilde{U}_{k+1}\tilde{N}_{k+}, \qquad A^T\tilde{U}_{k}=\tilde{V}_{k}\tilde{T}_{k},
\end{equation}
where $\tilde{T}_k \in \R^{k \times k}$ is upper triangular matrix with the entries $\tilde{t}_{i,j}$, 
$\tilde{N}_{k+} \in \R^{(k+1) \times k}$ is upper Hessenberg with the entries $\tilde{n}_{i,j}$, and 
$\tilde{Z}_k=[M_1^{-1}\tilde{v}_1,\ldots, M_k^{-1}\tilde{v}_k]$. 
The approximate solution $\tilde{x}_k$ of FLSQR is then given by the vectors $\tilde{z}_i, i=1, \dots, k$, as
\begin{equation}\label{FGKsolution}
\tilde{x}_k = x_0 + \tilde{Z}_k \tilde{y}_k = \argmin_{x \in x_0 + \textit{span}\{\tilde{z}_1, \dots, \tilde{z}_k \}}{\norm{b-Ax}},
\end{equation} 
where $\tilde{y}_k$ is the least squares solution of $\tilde{N}_{k+} \tilde{y}_k \approx \beta_1e_1$. 

FCGLS in Algorithm \ref{AlgorithmFCGLS} is adopted from \cite{nonNegative}, with a full derivation provided in the original paper. Due to the iteration-dependent preconditioner, a full recurrence on the update (direction) vectors $p_k$  is necessary to preserve their $A^TA$-orthogonality, see line $15$. The other computations remain in $3$-term recurrences. 
As a result, FCGLS is computationally less demanding than the FLSQR algorithm based on FGK. 
Note that by \cite[Equation (1)]{nonNegative}, FCGLS searches the approximate solution of~\eqref{theProblem} in the subspace generated by the vectors $\hat{s}_i, i=0, \dots, k-1$, such that
\begin{equation}\label{FCGLSsolution}
\hat{x}_k = \argmin_{x \in x_0 + \textit{span}\{\hat{s}_0, \dots, \hat{s}_{k-1}\} }{\norm{b-Ax}}.
\end{equation} 

\begin{minipage}[t]{0.46\textwidth}
\begin{algorithm}[H]
	\caption{Flexible GK \cite{flexibleKrylovLp}}
	\label{AlgorithmFGK}
	\begin{algorithmic}[1]
		\STATE{Input: }{$A,b,x_0, M_1, M_2, \ldots$}\\
		\STATE{$r_0=b-Ax_0$}\\
		\STATE{${\beta}_1=\norm{r_0}$}\\
		\STATE{$\tilde{u}_1=r_0/{\beta}_1$}\\	
		\FOR{$k =1, 2, \ldots$}
		\STATE{$\tilde{v}_k=A^T\tilde{u}_k$}\\
		\FOR{$j =1, 2, \ldots, k-1$}
		\STATE{$\tilde{t}_{j,k}=\tilde{v}_k^T\tilde{v}_j$}
		\ENDFOR\\
		\STATE{$\tilde{v}_k=\tilde{v}_k-\sum_{j=1}^{k-1}\tilde{t}_{j,k}\tilde{v}_j$}\\
		\STATE{$\tilde{t}_{k,k}=\norm{\tilde{v}_k}, \ \tilde{v}_k=\tilde{v}_k/\tilde{t}_{k,k}$}\\
		\STATE{	$\tilde{z}_k=M_k^{-1}\tilde{v}_k$}\\
		\STATE{$\tilde{u}_{k+1}=A\tilde{z}_k$}\\
		\FOR{$j=1, 2, \ldots, k$}
		\STATE{$\tilde{n}_{j,k}=\tilde{u}_{k+1}^T\tilde{u}_j$}\\
		\ENDFOR \\
		\STATE{$\tilde{u}_{k+1}=\tilde{u}_{k+1}-\sum_{j=1}^{k}\tilde{n}_{j,k}\tilde{u}_j$}\\
		\STATE{$\tilde{n}_{k+1,k}=\norm{\tilde{u}_{k+1}}$}\\
		\STATE{$\tilde{u}_{k+1}=\tilde{u}_{k+1}/\tilde{n}_{k+1,k}$}	
		\ENDFOR				
	\end{algorithmic}
\end{algorithm}
\end{minipage}
\hfill
\begin{minipage}[t]{0.46\textwidth}
	\begin{algorithm}[H]
		\caption{Flexible CGLS \cite{nonNegative}}
		\label{AlgorithmFCGLS}
		\begin{algorithmic}[1]
			\STATE{Input: }{$A,b,x_0, M_1,M_2, \ldots$}\\	
			\STATE{$r_0=b-A x_0$}\\
			\STATE{	$s_0=A^Tr_0$}\\
			\STATE{$\hat{s}_0=M_1^{-1}s_0$}\\	
			\STATE{	$p_0=\hat{s}_0,\ q_0=Ap_0$}\\
			\FOR{$k =1,2,  \ldots$}
			\STATE{$\gamma_{k-1}=(r_{k-1},q_{k-1})/\norm{q_{k-1}}^2$}\\
			\STATE{	$\hat{x}_{k}=\hat{x}_{k-1}+\gamma_{k-1} p_{k-1}$}\\
			\STATE{	$r_{k}=r_{k-1}-\gamma_{k-1} q_{k-1}$}\\
			\STATE{	$s_{k}=A^Tr_{k}$}\\
			\STATE{$\hat{s}_{k}=M_{k+1}^{-1}s_k$}\\
			\FOR{$j =0, 1, \ldots, k-1$}
			\STATE{$\theta_{j}^{k-1}=-(A\hat{s}_{k},q_j)/\norm{q_j}^2$}\\
			\ENDFOR
			\STATE{$p_{k}=\hat{s}_{k}+\sum_{j=0}^{k-1} \theta_{j}^{k-1} p_j$}\\
			\STATE{$q_{k}=A\hat{s}_{k}+\sum_{j=0}^{k-1} \theta_{j}^{k-1} q_j$}\label{AlgorithmFCGLSLineAp}
			\ENDFOR
		\end{algorithmic}
	\end{algorithm}
\vspace{0.1mm}
\end{minipage}

\section{Fast Flexible LSQR}\label{Section3}
To derive a computationally efficient LSQR with flexible preconditioning, we first introduce a novel flexible modification of the Golub-Kahan bidiagonalization. It is based on fac\-tori\-zation-free preconditioned LSQR for a fixed preconditioner presented in \cite{ArridgePreconditionedLSQR}, which we now briefly describe. This method was derived from the right-preconditioned LSQR, but does not require splitting (factorization) of the preconditioner $M$. Assuming $M$ is symmetric positive definite, the key idea is the application of $M^{-1}$-weighted inner products. The core of the bidiagonalization recurrence is, for $k=1, 2, \dots$,
\begin{align}
\breve{\beta}_{k}&\breve{u}_{k}=A\breve{z}_{k-1}-\breve{\alpha}_{k-1}\breve{u}_{k-1} \label{reccurenceUArridge}\\
\breve{v}_k&=A^T\breve{u}_{k}-\breve{\beta}_{k}\breve{v}_{k-1}, \label{reccurenceVArridge}
\end{align}
followed by computing $ \breve{z}_{k}=M^{-1}\breve{v}_{k}$ and scaling with the weighted inner product $\breve{\alpha}_{k}=(\breve{v}_{k},\breve{v}_{k})_{M^{-1}}^{\frac12}=(\breve{v}_{k},\breve{z}_{k})^{\frac12},\  \breve{v}_{k}=\breve{v}_{k}/\breve{\alpha}_{k},\ \breve{z}_{k}=\breve{z}_{k}/\breve{\alpha}_{k}$.
 The algorithm generates matrices 
$\breve{Z}_k$ and $\breve{U}_k$ with columns forming an $M$-orthonormal basis of $\mathcal{K}_k(M^{-1}A^TA,M^{-1}A^Tr_0)$ and
orthonormal basis of $\mathcal{K}_k(AM^{-1}A^T,r_0)$, respectively. Further, it computes the matrix $\breve{V}_k$ having 
$M^{-1}-$orthonormal columns such that
\begin{equation}\label{splitLSQR}	
  A\breve{Z}_k=\breve{U}_{k+1}\breve{L}_{k+}, \qquad A^T\breve{U}_k=\breve{V}_k\breve{L}_k^T,
\end{equation}
where $\breve{Z}_k = M^{-1}\breve{V}_k$ and $\breve{L}_k \in \R^{k \times k}, \breve{L}_{k+} \in \R^{(k+1) \times k}$ are lower bidiagonal matrices. See \cite{ArridgePreconditionedLSQR, EvickaIvetkaEnumath} for detailed derivations. 
The $k$th approximate solution of factorization-free preconditioned LSQR then satisfies
$$
\breve{x}_k \in\ x_0 + span\{ \breve{z}_1, \dots, \breve{z}_k \}= x_0+ \mathcal{K}_k(M^{-1}A^TA,M^{-1}A^Tr_0). 
$$	

\subsection{Fast Flexible Golub-Kahan recurrence}\label{SectionFlexibleRecurrence}
Consider now a sequence of changing symmetric positive definite preconditioners $M_1, M_2, \ldots, M_k$. 
In order to ensure orthogonality of basis vectors, we replace short recurrence generating $\breve{u}_k$ in \eqref{reccurenceUArridge} by a full recurrence. This will be later essential for formulation of the projected problem, see equation~\eqref{projectedProblemFaFLSQR}. Short recurrence for $\breve{v}_k$ from \eqref{reccurenceVArridge} is kept, but the fixed preconditioner $M$ is replaced by the flexible one. This gives our novel Fast Flexible Golub-Kahan algorithm presented in Algorithm \ref{AlgorithmFastFlexibleRecurrence}. 

\begin{algorithm}[ht]
		\caption{Fast Flexible Golub-Kahan recurrence}
	\label{AlgorithmFastFlexibleRecurrence}
	\begin{algorithmic}[1]
	\STATE{Input: }{$A,b,x_0, M_1, M_2, \ldots $}\\
	\STATE{$r_0=b-Ax_0$}\\
	\STATE{$u_1=r_0$}\\ 
	\STATE{$\beta_1=\norm{u_1},\ u_1=u_1/\beta_1$}\\
	\STATE{$v_0=0$}\\
	\FOR{$k = 1, 2, \ldots$}
		\STATE{$v_{k}=A^T u_{k}-\beta_{k} v_{k-1}$}\\
		\STATE{$z_{k}=M^{-1}_{k} v_{k}$}\\
		\STATE{$\alpha_{k}=(z_{k},v_{k})^{1/2}$}\\
		\STATE{$z_{k}=z_{k}/\alpha_{k},\ v_{k}=v_{k}/\alpha_{k}$}\\			
		\STATE{$u_{k+1}=A z_k $}\\
		\FOR{$j = 1, \ldots, k$}
			\STATE{$n_{j,k}=u_{k+1}^T u_j$}\\
		\ENDFOR
		\STATE{$u_{k+1}=u_{k+1}-\sum_{j=1}^{k} n_{j,k}u_j$}	\\		
		\STATE{$\beta_{k+1}=\norm{u_{k+1}},\ u_{k+1}=u_{k+1}/\beta_{k+1}$}
	\ENDFOR
\end{algorithmic}
\end{algorithm}

Denote $V_k=[v_1, v_2, \ldots, v_k], \ U_k=[u_1, u_2, \ldots,u_k]$ and $Z_k=[z_1, z_2, \ldots,z_k]=[M^{-1}_1v_1,\linebreak M^{-1}_2v_2, \ldots, M^{-1}_kv_k]$ matrices of vectors generated by Algorithm~\ref{AlgorithmFastFlexibleRecurrence}. Consider matrices of coefficients $L_k$ and $N_{k+}$ in the form
\begin{equation}\label{coefficientMatrices}
	L_k=
	\begin{bmatrix}
		\alpha_1    &          &       &           \\
		\beta_2     & \alpha_2 &       &           \\
		&          \ddots & \ddots  &		\\
		&          &\beta_{k} & \alpha_k
	\end{bmatrix} \in \R^{k\times k}, \ 
		N_{k+}=
	\begin{bmatrix}
		n_{1,1}   &  n_{1,2}  & \ldots      & n_{1,k}          \\
		\beta_2     & n_{2,2} &    \ldots   &         n_{2,k}  \\
		&          \ddots & \ddots  &		\\
		&          &\beta_{k} & n_{k,k}\\		
		&          & & \beta_{k+1}
	\end{bmatrix}\in \R^{(k+1) \times k}.
\end{equation}
Here $\alpha_j > 0$ since $M_j^{-1}$ is symmetric positive-definite and $\beta_{j+1} >0$ due to normalization, $j=1, \dots, k$. 
Then directly
\begin{equation}\label{FLSQRdecomposition}	
		AZ_k=U_{k+1}N_{k+}, \qquad A^TU_k=V_kL_k^T,
\end{equation}
compare with \eqref{FGKdecomposition} and \eqref{splitLSQR}.
The relations \eqref{FLSQRdecomposition} are similar to \eqref{splitLSQR}, except that we have the upper Hessenberg matrix 
$N_{k+}$ instead of the
lower bidiagonal $\breve{L}_{k+}$. The matrix $U_k$ has orthonormal columns thanks to the long recurrence orthogonalization.
\begin{remark}
Note that if $M_i=M$ for $i=1,2,\ldots k$, then Algorithm \ref{AlgorithmFastFlexibleRecurrence} reduces to the factorization-free preconditioned LSQR \cite{ArridgePreconditionedLSQR}. In such a case, $N_{k+}=\breve{L}_{k+}$ and the columns of $Z_k = \breve{Z}_k$ 
are $M$-orthogonal.
\end{remark}

\subsection{Orthogonality properties}\label{SectionOrthogonality}
Despite the fact that the vectors $v_k$ are computed by short recurrences, the newly computed $v_k$ is orthogonal to $z_{k-1}=M_{k-1}^{-1}v_{k-1}$. Thanks to the orthogonality of the vectors $u_i$ it is also orthogonal to all the previous $z_i=M_i^{-1}v_i$,
$i=1,2,\ldots k-2$. This important property requires a proof.

\begin{theorem}\label{TheoremOrthogonalityV}
The vectors $v_k,\ k=2,3, \dots$, generated by the Algorithm \ref{AlgorithmFastFlexibleRecurrence} have the following orthogonality property
\begin{equation}\label{orthogonalityFLSQR}
 v_k \perp_{M_i^{-1}} v_i, \quad \text{ or equivalently } \quad	v_k \perp z_i,\quad \ i=1,2,\ldots, k-1.
\end{equation}
Consequently, the matrix $V_k^TZ_k$ is upper triangular.
	\begin{proof}	
		Let $k,i \in \N, \  1\leq i<k$, then from Algorithm \ref{AlgorithmFastFlexibleRecurrence}
		\begin{align*}
			(v_k,z_{i})=&(v_k,M_{i}^{-1}v_{i})=(\alpha_k(A^Tu_k-\beta_kv_{k-1}),M_{i}^{-1}v_{i})\\
			&=\alpha_k(A^Tu_k,z_{i})-\alpha_k\beta_k(v_{k-1},M_{i}^{-1}v_{i}) \\
			&=\alpha_k(u_k,Az_{i})-\alpha_k\beta_k(v_{k-1},M_{i}^{-1}v_{i})\\
			&=\alpha_k\left(u_k,\beta_{i+1}u_{i+1}+\sum_{j=1}^{i}n_{j,i}u_j\right)-\alpha_k\beta_k(v_{k-1},M_{i}^{-1}v_{i}).
		\end{align*}
		For $i=(k-1)$, we have
		\begin{align*}
			\left(u_k,\beta_ku_k+\sum_{j=1}^{k-1}n_{j,k-1}u_j\right)&=\beta_k(u_k,u_k)=\beta_k, \\
			(v_{k-1},M_{k-1}^{-1}v_{k-1}) &= 1,			
		\end{align*}
		thanks to the orthogonality and normalization of $u_j$ in lines 15-16 and scaling of $v_k$ in lines 9-10. 
    For $i<(k-1)$, orthogonality of $u_j$ gives
		\begin{align*}
			\left(u_k,\beta_{i+1}u_{i+1}+\sum_{j=1}^{i}n_{j,i}u_j\right)=0.
		\end{align*}
		In summary, for $k=2,3,\ldots$
		\begin{align}
			(v_k,z_{k-1})&=\alpha_k\beta_k-\alpha_k\beta_k=0, \label{equation_vkzk1} \\
			(v_k,z_{i})&=-\alpha_k\beta_k(v_{k-1},z_{i}), \quad \quad \quad \ i=1,2,\ldots, k-2. \label{equation_vkzi}
		\end{align}			
		We now proceed by induction. For $k=2$, $v_2\perp z_1$ directly from \eqref{equation_vkzk1}.
		Assume that $v_{k-1}\perp z_i$ for $i=1,2,\ldots, k-2$. Then \eqref{equation_vkzk1} and \eqref{equation_vkzi} imply
		$v_{k}\perp z_i,\ i=1,2,\ldots, k-1$.	Then immediately $V_k^TZ_k$ is upper triangular.
	\end{proof}	
\end{theorem}

\begin{figure}[ht]
	\centering	
	\includegraphics[width=6cm]{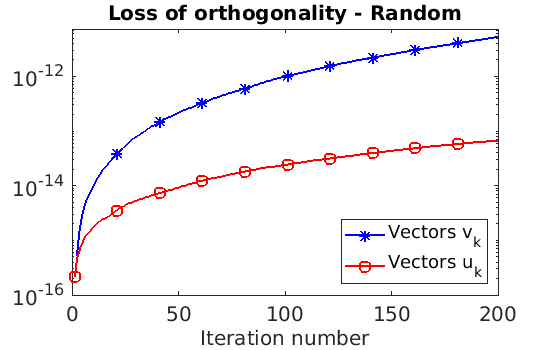}
	\includegraphics[width=6cm]{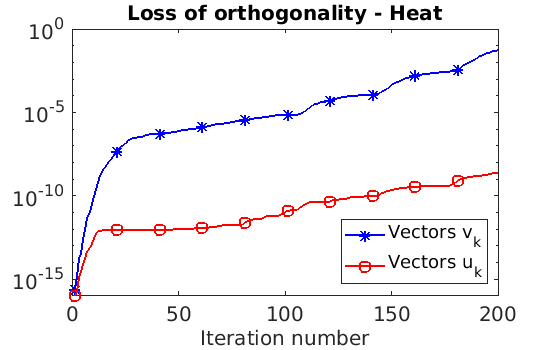}		
	\caption{Illustration of the loss of orthogonality in the Fast Flexible Golub-Kahan bidiagonalization. Left: a randomly generated 
	matrix $A\in \R^{n\times n}$ with $n=5000$, randomly generated $x^{exact}\in \R^{n}$, $b=A*x^{exact}+e$, where $e$ represents white noise with the noise level $\eta=1e-{4}$. Right: the \textit{heat} test problem described in section~\ref{Experiment1}. Both problems apply the same preconditioner described at the beginning of section~\ref{SectionExperiments}.}
	\label{FigureOGLoss}
\end{figure}

The orthogonality properties are illustrated on two examples in Figure~\ref{FigureOGLoss}.
We run $200$ iterations of the Fast Flexible Golub-Kahan bidiagonalization.
Since the vectors $u_k$ are 
in exact arithmetic orthonormal, the considered measure of orthogonality in finite precision 
is the classical $\|I - U_k^T U_k\|$, $k=1, \dots, 200$.
The matrix $V_k^TZ_k$ is in exact arithmetic upper triangular with $1$s on the main diagonal due to the
selected scaling on line $10$ of Algorithm \ref{AlgorithmFastFlexibleRecurrence}. Denote by $tril(V_k^TZ_k)$
the lower triangular part of the matrix including its diagonal. The orthogonality of the vectors $v_k$ is measured as
$\|I - tril(V_k^TZ_k)^T tril(V_k^TZ_k)\|$, $k=1, \dots, 200$. For the \textit{random} well posed problem (left) the orthogonality
is maintained well through iterations for both sets of vectors. For the ill-posed \textit{heat} test problem (right), 
the orthogonality is kept on reasonable level for a significant number of iterations.

\subsection{Derivation of Fast Flexible LSQR}\label{SectionFaFLSQR}
The vectors $z_1, z_2, \ldots,z_k$ generated by the Fast Flexible bidiagonalization provide a basis appropriate
for approximation of the solution of \eqref{theProblem} as we demonstrate later in section~\ref{SectionExperiments}.
Define the approximate solution $x_k$ as
\begin{equation}\label{solutionFaFLSQR}
	x_k = x_0+Z_ky_k = \argmin_{x \in x_0 + \textit{span}\{z_1, z_2, \ldots, z_k\} }{\norm{b-Ax}},
\end{equation}
compare with \eqref{FGKsolution}. Using \eqref{FLSQRdecomposition} and orthonormality of the columns of $U_k$, 
$$
\norm{b-Ax} = \norm{AZ_ky - r_0}=\norm{U_{k+1}N_{k+}y-U_{k+1}\beta_1e_1}=\norm{N_{k+}y-\beta_1e_1}.
$$
Thus the vector of coefficients is given by
$$
y_k = \argmin_{y\in \R^k}\norm{N_{k+}y-\beta_1e_1}.
$$
The approximate solution $x_k$ can be in each iteration $k$ computed by solving the small projected problem 
\begin{equation}\label{projectedProblemFaFLSQR}
N_{k+}y \approx \beta_1e_1
\end{equation}
in the least squares sense, followed by the transformation \eqref{solutionFaFLSQR}. 

Since $N_{k+}$ is upper-Hessenberg, \eqref{projectedProblemFaFLSQR} can be solved efficiently by the 
QR decomposition implemented through Givens rotations. According to \eqref{coefficientMatrices}, 
$N_{k+}$ is obtained by extending $N_{(k-1)+}$ by a zero row and one additional column. 
Thus an update formula for $x_k$ can be derived similarly as in GMRES \cite{GMRES}, see \cite{saad} for details.
For $k =1, 2, \dots$, consider the decomposition 
\begin{equation}\label{qrDecompositionL}
	Q_k [ N_{k+}, \beta_1 e_1 ]= 
	\begin{bmatrix}
		G_k & f_k \\
		& {\Phi}_{k+1}
	\end{bmatrix},
\end{equation}
where $Q_k$ is an orthonormal matrix, $G_k$ is upper triangular and $f_k$ is a column vector. Denote $D_k \equiv Z_kG_k^{-1} = [d_1, d_2, \ldots, d_k]$. 
Since $G_{k-1}$ is a submatrix of $G_{k}$, clearly $D_{k}=[D_{k-1},d_{k}]$, i.e.,
in the $k-$th iteration only the last column of $D_k$ has to be determined.
The update formula for the solution is then
\[ 
x_k = x_0+Z_kG_k^{-1}f_k = x_0+D_kf_k = x_{k-1}+f_k(k)d_k,
\]
where $f_k(k)$ is the $k-$th element of $f_k$. Then $G_k^TD_k^T=Z_k^T$ gives an expression 
for the update vector $d_k$, 
\begin{equation}\label{updateExpression}
	d_k=\frac{1}{g_{k,k}}\left( z_k - \sum_{i=1}^{k-1}g_{i,k}d_i \right),
\end{equation}
where $g_{i,k}$ is the element of $G_k$ at the position $(i, k)$.
In order to compute $d_k$, it is thus necessary to store all previous update 
vectors $d_i$, but we no longer store the vectors $z_i$ as for \eqref{solutionFaFLSQR}.
The whole algorithm is summarized in Algorithm \ref{AlgorithmFaFLSQR}.

\begin{algorithm}[ht]
	\caption{Fast Flexible LSQR}
	\label{AlgorithmFaFLSQR}
	\begin{algorithmic}[1]	
	\STATE{Input: }{$A,b,x_0, M_1, M_2, \ldots $}\\
	\STATE{$r_0=b-Ax_0$}\\
	\STATE{$u_1=r_0$}\\ 
	\STATE{$\beta_1=\norm{u_1},\ u_1=u_1/\beta_1$}\\
	\STATE{$v_0=0$}\\
	\FOR{$k = 1, 2, \ldots$}
		\STATE{$v_{k}=A^T u_{k}-\beta_{k} v_{k-1}$}\\
		\STATE{$z_{k}=M^{-1}_{k} v_{k}$}\\
		\STATE{$\alpha_{k}=(z_{k},v_{k})^{1/2}$}\\
		\STATE{$z_{k}=z_{k}/\alpha_{k},\ v_{k}=v_{k}/\alpha_{k}$}\\			
		\STATE{$u_{k+1}=A z_k $}\\
		\FOR{$j = 1, \ldots, k$} 
			\STATE{$n_{j,k}=u_{k+1}^T u_j$}
		\ENDFOR
		\STATE{$u_{k+1}=u_{k+1}-\sum_{j=1}^{k} n_{j,k}u_j $} 	\\		
		\STATE{$\beta_{k+1}=\norm{u_{k+1}},\ u_{k+1}=u_{k+1}/\beta_{k+1}$}\\
		\STATE{\mbox{update the QR decomposition in \eqref{qrDecompositionL}}  }\\
		\STATE{$d_{k}=\frac{1}{g_{k,k}}\left( z_{k} - \sum_{i=1}^{k-1}g_{i,k}d_i \right)$}	\\	
		\STATE{$x_{k}=x_{k-1}+f_k(k)d_{k}$}\\	
	\ENDFOR	
\end{algorithmic}
\end{algorithm}

\subsection{Hybrid FaFLSQR}

When solving \eqref{theProblem}, the ill-posedness of the problem and the ill-conditioning of the model matrix make the solution highly sensitive to perturbations in the data. Projection onto a Krylov subspace, such as through the Golub-Kahan bidiagonalization, acts as a form of regularization. However, the projected problem tends to inherit some of the original ill-posedness, and noise gradually propagates into the projected system; see e.g. \cite{NoiseHnetynkova}. Consequently, these methods frequently exhibit semiconvergence. This phenomenon is well known for standard Krylov subspace methods and it occurs also with flexible preconditioning, see section~\ref{SectionExperiments} for illustration. Thus, applying additional regularization to the projected problem is often advantageous; see \cite{Hansen} and the references therein. 

The minimization problem in FaFLSQR is similar to the one in FLSQR, compare the equations \eqref{solutionFaFLSQR} and \eqref{FGKsolution}. The derivation of Hybrid FaFLSQR therefore closely parallels that of Hybrid FLSQR~\cite{flexibleKrylovLp}. Shortly, applying Tikhonov regularization to the projected problem in FaFLSQR with a regularization parameter $\lambda_k \in \mathbb{R}$ gives Hybrid FaFLSQR. Mathematically, the obtained minimization in iteration $k$ is
\begin{equation}\label{hybridFLSQR} 
	\min_{y\in \R^k} \{ \norm{N_{k+} y - \beta_1 e_1}^2 + \lambda_k^2 \norm{y}^2 \}. 
\end{equation}
This could be further expanded to a more general form,
\begin{equation}\label{generalhybridFLSQR}
	\min_{y\in \R^k} \{ \norm{N_{k+} y - \beta_1 e_1}^2 + \lambda_k^2 \norm{J_k y}^2 \},
 \end{equation} 
where $J_k$ is a suitably chosen matrix, typically such that it enforces some conditions on the solution. For further details, see, e.g. \cite{flexibleKrylovLp}. In either case, the approximate solution is then recovered via $x_k = x_0 + Z_k y_k$. 

In section~\ref{SectionExperiments}, we demonstrate the behavior of Hybrid FaFLSQR through numerical experiments 
and compare it to Hybrid FLSQR.

\section{Comparison to other flexible algorithms}
When flexible preconditioning is applied in Krylov subspace methods, it is difficult to analytically express the approximation subspaces as they are no longer standard Krylov subspaces. Therefore, bases generated by various flexible methods are studied in this section. We show the difference between FaFLSQR and FLSQR. We prove that FaFLSQR is mathematically equivalent to FCGLS in the sense that they construct the same approximation in each iteration $k$.
Furthermore, computational cost of algorithms is compared.

\subsection{Comparison to FLSQR}
Recall that quantities generated by FLSQR are denoted by tilde. Consider the same initial approximation 
$x_0$ and preconditioners $M_1, \ldots, M_k$. FaLSQR (Algorithm~\ref{AlgorithmFastFlexibleRecurrence}) and FLSQR 
(Algorithm~\ref{AlgorithmFGK}) both reduce the original problem to the projected problem 
with an upper-Hessenberg matrix $N_{k+}$ and 
$\tilde{N}_{k+}$, respectively. The approximate solutions are linear combinations of basis vectors $z_k$ and $\tilde{z}_k$, respectively,
see \eqref{solutionFaFLSQR} and \eqref{FGKsolution}. 

Clearly, $\tilde{u}_1 = u_1$, $\textit{span}\{\tilde{v}_1\}= \textit{span}\{v_1\}$ (note the different scaling), and thus 
\begin{equation}\label{equations_z1}
	\textit{span}\{\tilde{z}_1\}= \textit{span}\{z_1\}=\textit{span}\{M_1^{-1}v_1\} .
\end{equation}
For $k>1$,
\begin{equation}\label{equations_zk}
\begin{aligned}
	z_k =&M_k^{-1}v_k=\frac{1}{\alpha_k}M_k^{-1}\left[ A^Tu_k-\beta_kv_{k-1} \right],
	\\
	\tilde{z}_k=&M_k^{-1}\tilde{v}_k=\frac{1}{t_{k,k}}M_k^{-1}\left[ A^T\tilde{u}_{k}-\sum_{i=1}^{k-1} \tilde{t}_{i,k}\tilde{v}_i\right].\\
\end{aligned}
\end{equation}
Put $k = 2$. From \eqref{equations_z1}, we get $\tilde{u}_2 = u_2$ and 
\begin{equation}\label{equations_z2}
	\begin{aligned}
		z_2 &= \frac{1}{\alpha_2}M_2^{-1}\left[ A^Tu_2 - \beta_2 v_1 \right], \quad 
		\tilde{z}_2 &= \frac{1}{\tilde{t}_{2,2}}M_2^{-1}\left[ A^Tu_2 - \tilde{t}_{1,2} \tilde{v}_1 \right].
	\end{aligned}
\end{equation}
Using \eqref{equations_z1} and \eqref{equations_z2}, if $\beta_2 v_1 = \tilde{t}_{1,2} \tilde{v}_1$ or $M_2^{-1}v_1$ is a scalar multiple 
of $M_1^{-1}v_1$, then $\textit{span}\{z_1, z_2\} = \textit{span}\{\tilde{z}_1, \tilde{z}_2\}$.
Recall that $\tilde{v}_2$ is orthogonal to $\tilde{v}_1$, whereas $v_2$ is $M_2^{-1}$-orthogonal to $v_1$. Therefore, the condition $\beta_2 v_1 = \tilde{t}_{1,2} \tilde{v}_1$ is satisfied when $v_2$ is both orthogonal and $M_1^{-1}$-orthogonal to $v_1$. This happens, e.g., when 
$M_1^{-1}$ is a scalar multiple of the identity matrix. 

\begin{figure}[ht]
	\centering	
	\includegraphics[width=13cm]{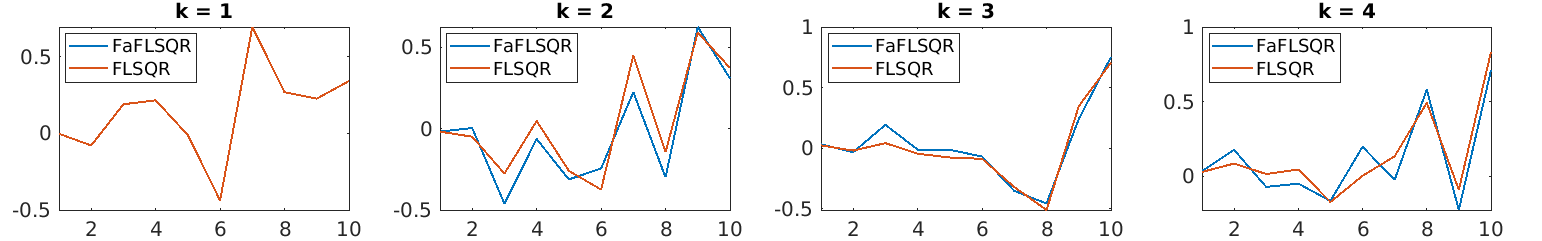}	
	\caption{Comparison of normalized basis vectors $z_k$, $\tilde{z}_k$ generated by FaFLSQR and FLSQR, respectively. Here \textit{random}
	test problem has the size $n=10$.}
	\label{FigureBases}
\end{figure}

If neither of the two conditions above is satisfied, then $\textit{span}\{z_1, z_2\} \neq \textit{span}\{\tilde{z}_1, \tilde{z}_2\}$ and FLSQR and FaFLSQR can generate different approximations. To illustrate this, we provide a simple example. The setting is the same as for the \textit{random} test problem in Figure~\ref{FigureOGLoss}, except $n=10$ here to allow visualization. An iteration dependent preconditioner for both algorithms is $M_i^{-1}=\textit{diag}(i^2, (i+1)^2, \ldots, (i+9)^2)$ 
(for practical preconditioners see section 5). Figure~\ref{FigureBases} illustrates that even though the generated normalized 
vectors $z_k$, $\tilde{z}_k$ for $k>1$ are close, they are not equal.
Moreover, they generate different bases, see ranks of the concatenated matrix $[Z_k,\tilde{Z}_k]$ in Table~\ref{TableRank}. In the first iteration, $z_1$ and $\tilde{z}_1$ are collinear, thus the rank is $1$. Further, the rank increases by two in each iteration. Consequently, the computed approximations $x_k$, $\tilde{x}_k$ are different, see line $3$ of Table~\ref{TableRank}.
\begin{table}[htbp]
	\footnotesize
	\caption{Rank of the concatenated matrix $[Z_k,\tilde{Z}_k]$ and relative difference between the approximations $x_k$, $\tilde{x}_k$, for
	iterations $k=1, \dots, 5$, \textit{random} test problem with $n=10$.}
	\label{TableRank}
	\begin{center}
		\begin{tabular}{|c|c|c|c|c|c|c|} \hline
			Number of iterations (k)	&  
			1  &   2  &  3     &  4     &    5  \\ \hline
			rank($[Z_k,\tilde{Z}_k])$ &		1  &    3  &   5   &   7  &   9  \\
			$\norm{x_k-\tilde{x}_k}/\norm{x_k}$ &  $9.1490e-17$  &  0.0242 & 0.0389 &  0.0523  &  0.0494  \\ \hline
		\end{tabular}    
	\end{center}
\end{table}

\begin{figure}[ht]
	\centering	
	\includegraphics[width=6cm]{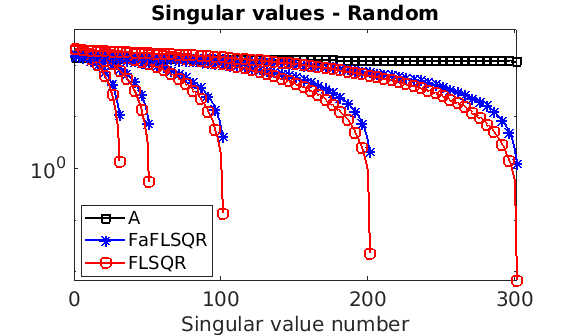}	
	\includegraphics[width=6cm]{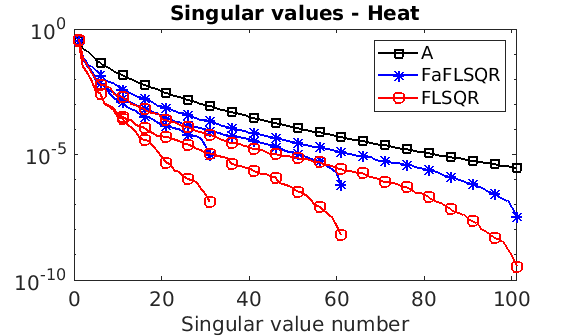}	
	\caption{Comparison of singular values of $A$, and their approximations given by the singular values of $\tilde{N}_{k+}$ from FLSQR and $N_{k+}$ from FaFLSQR. The test problems are the same as in Figure~\ref{FigureOGLoss}. Left: $k=30, 50, 100, 200, 300$. Right: $k=30, 60, 100$. Every fifth singular value is highlited by a marker.}
	\label{FigureSVD}
\end{figure}

Different properties of FaLSQR and FLSQR can be seen also when comparing singular values of $N_{k+}$ and $\tilde{N}_{k+}$. 
It is well known that for the basic Golub-Kahan bidiagonalization (underlying LSQR) the singular values of the projected matrix 
$\bar{L}_{k+}$ approximate the singular values of the original matrix $A$. This property does not directly carry over to the flexible preconditioned variants, but it is still interesting to study $N_{k+}$, $\tilde{N}_{k+}$ and $A$.
Figure~\ref{FigureSVD} depicts singular values for the test problems \textit{random} (left) 
and \textit{heat} (right) as in Figure~\ref{FigureOGLoss}, for selected iterations $k$. 
For the well-conditioned \textit{random} problem, the singular values of all three matrices are close. 
For the ill-conditioned \textit{heat} problem, the singular values of $N_{k+}$ (FaFLSQR) provide better approximations 
than the ones of $\tilde{N}_{k+}$ (FLSQR).

\subsection{Equivalency with FCGLS}\label{SectionEquivalency}
The equivalency of the standard CGLS and LSQR can be proven easily by showing that both methods minimize the same functional (residual norm) over the same Krylov subspace in each iteration $k$. Similarly, we prove equivalency of the novel FaFLSQR given in Algorithm \ref{AlgorithmFastFlexibleRecurrence} and FCGLS presented in \cite{nonNegative}, see also Algorithm \ref{AlgorithmFCGLS}. 

FaFLSQR searches the $k$th approximate solution by minimizing  \eqref{solutionFaFLSQR} and FCGLS by minimizing \eqref{FCGLSsolution}. Thus
it is enough to prove that $z_k$ in FaFLSQR is a scalar multiple of $\hat{s}_{k-1}$ in FCGLS, $k=1,2, \dots$. First, the following theorem shows that the vectors $s_k$ have similar orthogonality property as given in Theorem \ref{TheoremOrthogonalityV} for the vectors $v_k$.

\begin{theorem}
	The vectors $s_k, k=1,2, ...$, generated by the Algorithm \ref{AlgorithmFCGLS} \newline (FCGLS) have the following orthogonality property
	\[
	s_k \perp_{M^{-1}_{i+1}} s_i,\ i=0,1,\ldots, k-1.
	\]
	\begin{proof}
		See the Appendix.
	\end{proof}
\end{theorem}

\begin{theorem}\label{TheoremBasisEquality}
The vectors $z_k$ generated by the Algorithm \ref{AlgorithmFastFlexibleRecurrence} (FaFLSQR) are scalar multiples of 
the vectors $\hat{s}_{k-1}$ generated by the Algorithm \ref{AlgorithmFCGLS} (FCGLS), i.e.
	\[
	z_k=\eta_k \hat{s}_{k-1}
	\] 
	for some $\eta_k \in \R, \eta_k\neq0,\ k=1,2,\ldots $.
\begin{proof}
Since $z_k=M^{-1}_kv_k$ and $\hat{s}_{k-1}=M^{-1}_{k}s_{k-1}$, $k=1,2,\ldots$, it is enough to prove that $v_k$ are scalar multiples of $s_{k-1}$. From Algorithms \ref{AlgorithmFastFlexibleRecurrence} and \ref{AlgorithmFCGLS}, 
$$
v_1=\frac{1}{\alpha_1 \beta_1}A^Tr_0, \quad s_0=A^Tr_0, \quad \quad \mbox{giving} \quad \quad v_1= \frac{1}{\alpha_1\beta_1} \, s_0.
$$
We proceed by induction. Assume that $v_i=\eta_i s_{i-1}$ for $\eta_i \neq 0, \eta_i \in \R,$ for all $i=1,2,\ldots, k$. Then
		\begin{align*}		
			v_{k+1}&=\frac{1}{\alpha_{k+1}}\left( A^Tu_{k+1}-\beta_{k+1}v_{k} \right)\\
			&\quad \quad \quad=\frac{1}{\alpha_{k+1}}\left( \frac{1}{\beta_{k+1}}A^T\left(Az_{k}-\sum_{j=1}^{k}n_{j,k} u_j  \right) -\beta_{k+1}v_{k} \right)\\
			&\quad \quad \quad=\frac{1}{\alpha_{k+1}\beta_{k+1}}\left( A^TAz_{k}- \sum_{j=1}^{k}n_{j,k} (\alpha_j v_j+\beta_jv_{j-1}) -\beta_{k+1}^2v_{k} \right), \\		
			s_k&=A^Tr_k=A^T(r_{k-1}-\gamma_{k-1}q_{k-1})\\
			&\quad \quad \quad=s_{k-1}-\gamma_{k-1}A^Tq_{k-1}\\
			&\quad \quad \quad=s_{k-1}-\gamma_{k-1}A^T\left(A\hat{s}_{k-1}-\sum_{j=0}^{k-2}\theta^{k-2}_{j}q_j \right)\\
			&\quad \quad \quad=s_{k-1}-\gamma_{k-1}A^TA\hat{s}_{k-1}+\gamma_{k-1}\sum_{j=0}^{k-2}\theta^{k-2}_{j}A^Tq_j\\
			&\quad \quad \quad=s_{k-1}-\gamma_{k-1}A^TA\hat{s}_{k-1}+\gamma_{k-1}\sum_{j=0}^{k-2}\theta^{k-2}_{j}\frac{1}{\gamma_{j}}(s_{j}-s_{j+1}).
		\end{align*}
		Thus
		\begin{align}\label{equationSpanvksk}
			\begin{split}
				v_{k+1} &\in \textit{span}\{v_1,v_2, \ldots, v_{k}, A^TAM_k^{-1}v_k\},\\
				s_{k}   &\in \textit{span}\{s_0, s_1, \ldots, s_{k-1}, A^TAM_k^{-1}s_{k-1}\},
			\end{split}
		\end{align}
		where $\textit{span}\{ v_1, v_2, \ldots, v_{k+1} \}= \textit{span}\{ s_0, s_1, \ldots, s_{k} \}$ by the induction assumption. \newline    Combining this with the fact that $v_{k+1} \perp_{M_{i}^{-1}}v_i,\ i=1,\ldots, k$, and 
		$s_k \perp_{M_{i+1}^{-1}}s_i$, $i=0,\ldots, k-1$, yields
		\begin{equation}
			v_{k+1}=\eta_{k+1}s_k
		\end{equation}
		for some $\eta_{k+1} \neq 0, \ \eta_{k+1} \in \R$, which finishes the proof.
	\end{proof}
\end{theorem}
\begin{remark}
Note that the coefficients $\eta_k$ can be derived explicitly. Specifically, we have 
	$\eta_k=(-1)^{(k+1)}\frac{1}{\alpha_1\beta_1}\Pi_{i=2}^{k}\left(\frac{1}{\alpha_i\beta_i\gamma_{i-2}}\right)$. 
\end{remark}
Numerical illustration of the equivalency between FCGLS and FaFLSQR is given in section~\ref{Experiment1}.

\subsection{Computational cost}
Now we compare computational costs of FaFLSQR (Algorithm~\ref{AlgorithmFaFLSQR}), FLSQR (Algorithm~\ref{AlgorithmFGK}) and FCGLS (Algorithm~\ref{AlgorithmFCGLS}). Let $A\in\R^{n\times n}$ for simplicity of exposition, generalization to rectangular matrices is straightforward.

\paragraph*{Matrix-vector multiplications}
The potentially most computationally demanding operations at each iteration $k$ of all the algorithms are two matrix-vector multiplications with $A$ and $A^T$, respectively. In terms of the number of floating point operations these together cost $2\times O(n^2)$. However, if $A$ is sparse or, for example, represents a point spread function, the computational cost reduces significantly to $2\times O(n \log(n))$ or $2\times O(cn)$ (with the constant $c\in\R, c>0$, depending on the sparsity of $A$). The multiplication with the preconditioning matrix $M_k$ occurs ones in each iteration and its cost is negligible since $M_k$ is typically a diagonal matrix.
 
\paragraph*{Recurrences}
FaFLSQR and FCGLS compute one full recurrence on the vectors $u_k$ and $q_k$, respectively. In iteration $k$ this costs $O(kn)$. The rest of the computations relies on short recurrences, in particular determination of the vectors $v_k$ in FaFLSQR and $r_k$ in FCGLS. The cost of short recurrences is negligible compared to the cost of long recurrences. The FLSQR algorithm computes two long recurrences on the vectors $u_k$ and $v_k$, the cost in iteration $k$ is thus $2 \times O(kn)$.
 
\paragraph*{Solution update}
The approximate solution typically needs to be computed in each iteration since it is often used to determine $M^{-1}_{k+1}$. For FaFLSQR and FLSQR this can be done efficiently using the update formula \eqref{updateExpression}. In iteration $k$, this costs $O(kn)$ operations. Alternatively, the approximate solution ${x}_k$, $\tilde{x}_k$ can be computed by solving the system with the projected matrix $N_{k+}$ and $\tilde{N}_{k+}$, respectively, and then using \eqref{solutionFaFLSQR} with the cost $O(k^2)$ + $O(kn)$. The latter option is used also for the Hybrid versions of the algorithms. In FCGLS, the direction vectors $p_k$ are used to update the approximate solution $\hat{x}_k$ in each iteration. In iteration $k$, computing the vector $p_k$ costs $O(kn)$.
 
\paragraph*{Storage}
Both FaLSQR and FCGLS require storing two full sets of vectors. For FaFLSQR these are the vectors $u_k$ necessary to compute the long recurrence, and either the update vectors $d_k$ or the basis vectors $z_k$ in order to recover $x_k$, see derivations in section~\ref{SectionFaFLSQR}. For FCGLS, the stored vectors are $q_k$ and $p_k$.
FLSQR requires storing three full sets of vectors, particularly $\tilde{u}_k$ and $\tilde{v}_k$ for the two long recurrence orthogonalizations, and one set of vectors to recover the approximate solution $\tilde{x}_k$ (update vectors for the update formula, or the basis vectors $\tilde{z}_k$).
  
\paragraph*{Discussion}
The computational and storage costs of FaFLSQR and FCGLS are similar. On the contrary, FLSQR requires evaluation of one additional long recurrence and consequently one extra full set of stored vectors. The computational cost of this additional recurrence might be negligible compared to the cost of the matrix-vector multiplications with $A$ for small $k$. However, for structured or sparse matrices, the cost of multiplications reduces significantly, while the cost of the long recurrences remains the same. Moreover, matrix-vector multiplications are typically easy to paralelize unlike the computation of the dot products in the long recurrences. 
This implies that computational efficiency of FaFLSQR can be remarkably higher than of FLSQR for many inverse problems. \\

\begin{figure}[ht]
	\centering	
	\label{FigureDenseSparseTiming}
	\includegraphics[width=6cm]{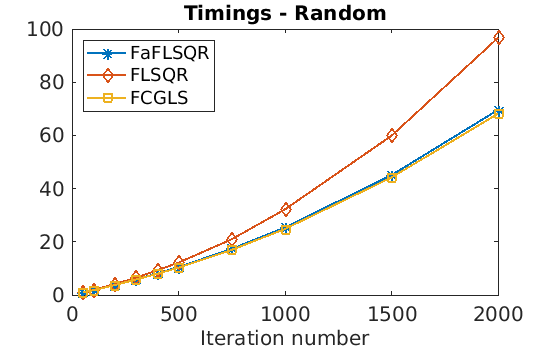}
	\includegraphics[width=6cm]{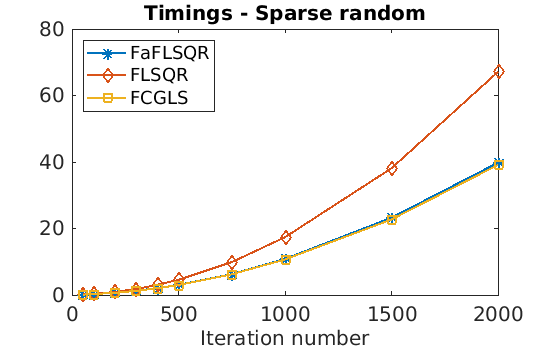}	
	\caption{Comparison of the computational time (measured in seconds) required by FaFLSQR, FLSQR and FCGLS, with respect 
	to the number of iterations $k$. Left: a system with a square random matrix $A\in\R^{5 \, 000 \times 5 \, 000}$ (the same as in Figure~\ref{FigureOGLoss}). Right: a system with a sparse random matrix $A\in\R^{5 \, 000 \times 5 \, 000}$ with the density of nonzero elements $5 \%$.}	
\end{figure}

We illustrate this in Figure~\ref{FigureDenseSparseTiming} and Table~\ref{TableDenseSparse} on two examples with random matrices -
one full and one sparse, $A\in\R^{5\, 000  \times 5 \, 000}$. The specification of the computer used for the experiments and implementation details are provided at the beginning of section~\ref{SectionExperiments}. In order to compare computational time $T$ taken by FaFLSQR to computational time $T_o$ of other algorithms, define the value $\tau$ representing the speed up obtained by FaFLSQR as
\begin{equation}\label{tauDefinition}
\tau= 100 *\left(1-\frac{T}{T_o}\right) \% .
\end{equation}
For a sparse matrix $A$ and a larger number of iterations, the difference in Figure~\ref{FigureDenseSparseTiming} (right) and Table~\ref{TableDenseSparse} (last row) is more pronounced. Even for a small number of iterations $k=50$ (compared to the size of $A$), FaFLSQR is about $15\%$ faster than FLSQR. Further results are presented in section~\ref{SectionExperiments}.

\begin{table}[htbp]
	\footnotesize
	\caption{Speed up $\tau$ of FaLSQR in comparison to FLSQR with respect to the number of iterations $k$. 
	The value $\tau$ is defined in \eqref{tauDefinition}. }\label{TableDenseSparse}
	\begin{center}
		\begin{tabular}{|c|c|c|c|c|c|c|c|c|c|} \hline
			Number of iters.	&  
			50  &   100  &  200     &  400     &    750  &   1000  & 1500   &  2000 \\ \hline
 			Dense m. &		4.9 \% &    5.5 \% &   7.0  \% &   13.2 \% &   17.4  \% & 21.0 \% &  24.7  \% & 28.0 \%  \\
			Sparse m. &  14.3 \% &   18.4 \% & 25.2 \% &   32.7 \% &   36.4  \% &    37.6   \% &   39.0   \% &   40.7 \%      \\ \hline
		\end{tabular}    
	\end{center}

\end{table}

\section{Numerical experiments}\label{SectionExperiments}
In this section, we numerically illustrate properties of FaFLSQR and compare it to FLSQR and FCGLS. Three different standard benchmark problems are used to provide comprehensive comparison. For simplicity, the initial approximation $x_0$ is a zero vector in all experiments. For flexible preconditioning, consider a commonly used diagonal preconditioner 
$M_1=I$, $M_k=diag(1/ \Psi(x_{k-1})), k>1$, where $\Psi=\max\{\abs{x_{k-1}}, tol\}$ with the safety parameter $tol=10^{-10}$, see also \cite{flexibleKrylovLp}.

We use our own implementation of all algorithms. In order to provide a fair comparison, both FaFLSQR and FLSQR include the update formula derived in section~\ref{SectionFaFLSQR}. All long recurrence orthogonalizations are done by the Modified Gram-Schmidt algorithm; cf. 
\cite{GramSchmidt} for history and variants. 
The implementations are avilable at: https://gitlab.com/havelkova-eva/fast-flexible-lsqr.
Experiments were performed on a laptop with a 11th Gen Intel(R) Core(TM) i7-1165G7 @ 2.80GHz processor and 16 GB DDR4 Synchronous 
3200 MHz Memory.

\subsection{Experiment 1}\label{Experiment1}
The first experiment illustrates numerically the equivalency of FaFLSQR and FCGLS proved in section~\ref{SectionEquivalency}. Consider the \textit{heat} test problem from~\cite{Regutools} with the dimension $n=500$. The right hand side $b$  is polluted by Gaussian white noise with the noise level $\eta=1e-{4}$. 
 
First, we compare the quality of the computed approximations by means of the relative true error $\norm{x_k-x^{exact}}/\norm{x^{exact}}$
for $k=1,2, \dots$. Recall that $x^{exact}$ is the exact solution of the unknown noise-free problem $A x^{exact} = b^{exact}$.
\begin{figure}[ht]
	\centering	
	\includegraphics[width=6cm]{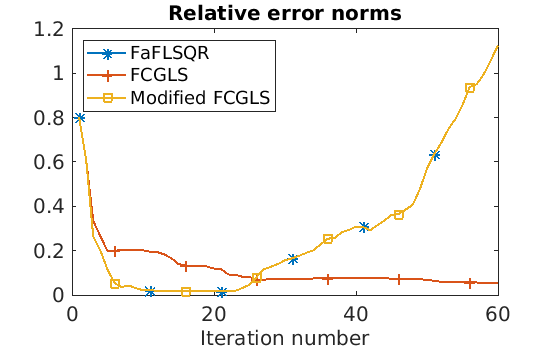}
	\includegraphics[width=6cm]{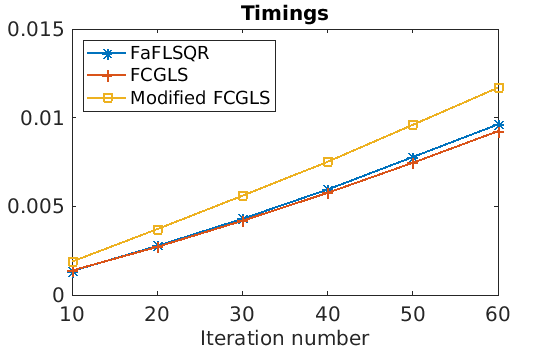}
	\label{FigureFCGLSError}
	\caption{Experiment 1. Relative error norms (left) and computational time measured in seconds (right) with respect to the number of      
   iterations $k$, for FaFLSQR, FCGLS and Modified FCGLS on the \textit{heat} test problem.}
\end{figure}
Figure~\ref{FigureFCGLSError} (left) shows error norms with respect to $k$. The observed significant delay of convergence in FCGLS is caused by rounding errors in finite precision computations. 
This phenomenon is well described for standard conjugate gradients (CG) \cite{CGLS} based on short recurrences; see \cite[Chapter 5]{Strakos}
for detailed discussion and further references. Therefore, in order to numerically illustrate the equivalency of FCGLS and FaFLSQR, consider also a modified implementation of FCGLS (denoted as Modified FCGLS) where the line \ref{AlgorithmFCGLSLineAp} of Algorithm \ref{AlgorithmFCGLS} is replaced by direct computation of $q_k=Ap_k$. This does not affect analytical properties of FCGLS, and thus Modified FCGLS is mathematically equivalent to FaFLSQR. In Figure~\ref{FigureFCGLSError} (left), the Modified FCGLS attains error norms similar to FLSQR. However, the modification introduces an additional matrix-vector product with $A$ and thus increases computational cost, see Figure~\ref{FigureFCGLSError} (right). 

Theorem~\ref{TheoremBasisEquality} states that the basis vectors $z_k$ from FaLSQR and $\hat{s}_{k-1}$ from FCGLS are equal 
up to scaling. Consequently, also the update vectors $d_k$ of the approximate solutions $x_k$ 
(see line 19 of Algorithm~\ref{AlgorithmFaFLSQR}) and $p_{k-1}$ of
$\hat{x}_k$ (see line 8 of Algorithm~\ref{AlgorithmFCGLS}), must be collinear in each iteration. 
\begin{figure}[ht]
	\centering	
	\includegraphics[width=13cm]{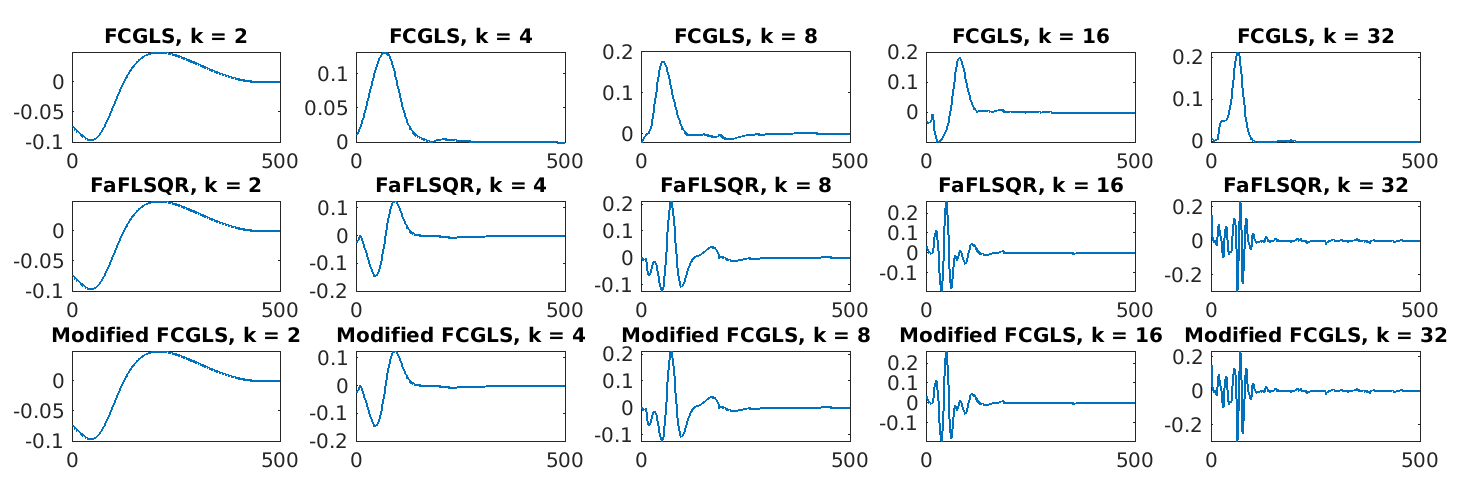}
	\label{FigureFCGLSBasesCompare}
	\caption{Experiment 1. Comparison of selected normalized basis vectors $z_k$ (FaFLSQR) and $\hat{s}_{k-1}$ (FCGLS and Modified FCGLS) for the \textit{heat} test problem.}
\end{figure}

\begin{figure}[ht]
	\centering	
	\includegraphics[width=13cm]{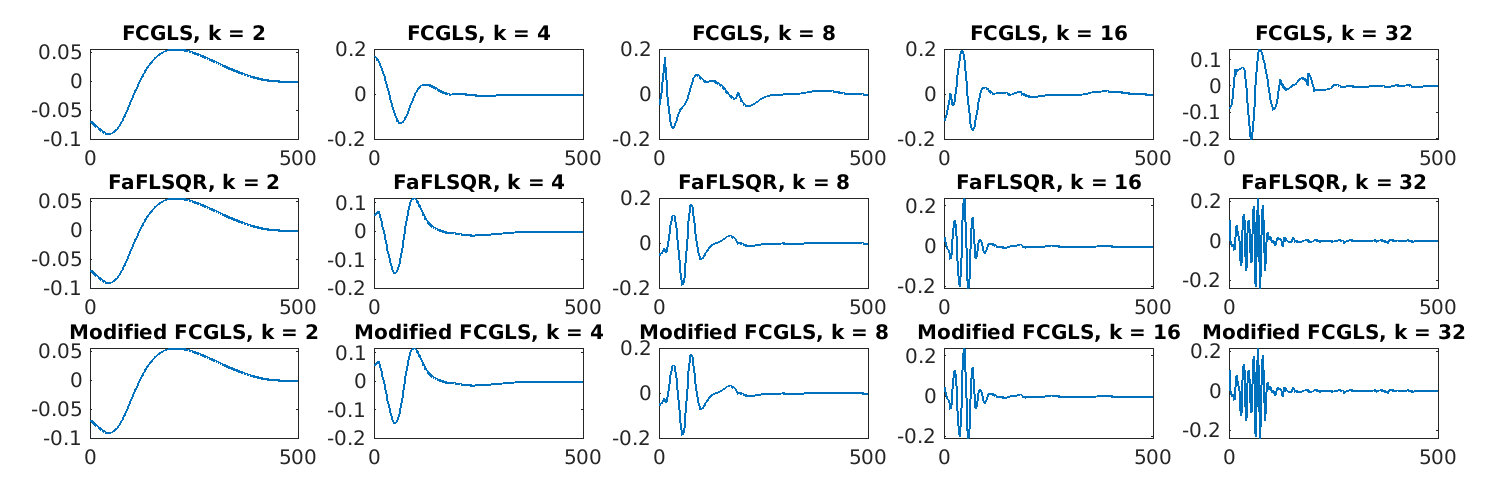}
	\label{FigureFCGLSUpdatesCompare}
	\caption{Experiment 1. Comparison of selected normalized update vectors $d_k$ (FaFLSQR) and $p_{k-1}$ (FCGLS and Modified FCGLS) for the \textit{heat} test problem.}
\end{figure}

Figure~\ref{FigureFCGLSBasesCompare} compares normalized vectors $z_k$ from FaFLSQR with normalized vectors $\hat{s}_{k-1}$ computed by FCGLS and Modified FCGLS, respectively. The vectors determined by FCGLS start to differ from $z_k$ significantly in iteration $k=4$ due to the delay of convergence. For Modified FCGLS, the bases vectors $\hat{s}_{k-1}$ correspond nicely to $z_k$. Similar behavior is present for the normalized update vectors $d_k$ of FaFLSQR and normalized update vectors $p_{k-1}$ of FCGLS and Modified FCGLS, see Figure~\ref{FigureFCGLSUpdatesCompare}. 

\subsection{Experiment 2}\label{Experiment2}
\begin{figure}[ht]
	\centering	
	\includegraphics[width=13cm]{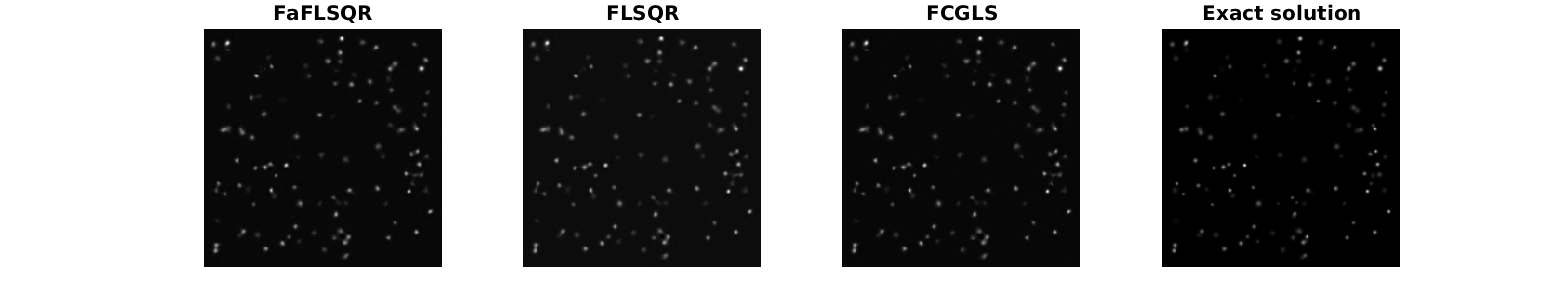}
	\includegraphics[width=13cm]{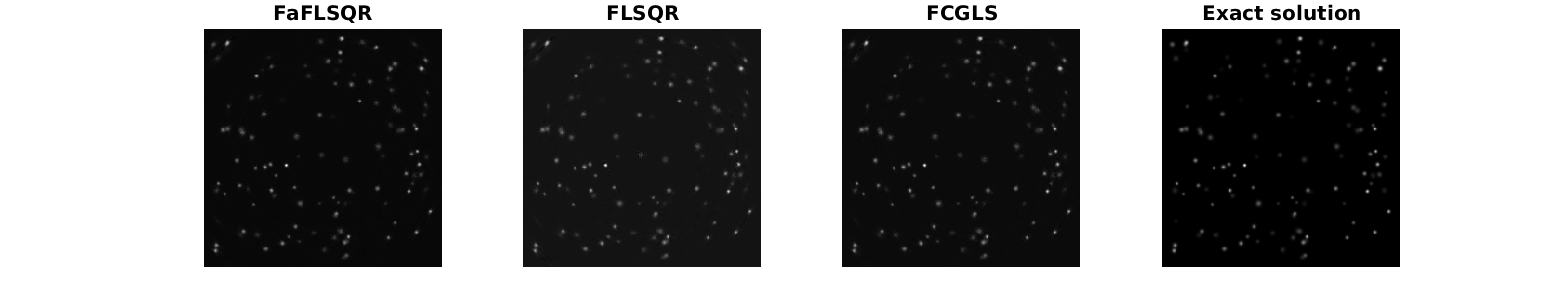}
	\label{FigureBlurSolutions}
	\caption{Experiment 2. Comparison of the exact solution and approximations with the lowest error norms obtained 
	by FaLSQR, FLSQR and FCGLS. Upper images correspond to the Gaussian blur, lower images to the Rotational blur in
	the test problem \textit{PRBlur}, $\eta=5e-{2}$.}
\end{figure}

In this experiment, we compare FCGLS and FLSQR to FaFLSQR. Consider the two dimensional image deblurring \textit{PRBlur} problem from IRtools toolbox \cite{IRTools} with two different point spread functions - spatially invariant Gaussian blur with periodic boundary conditions, and spatially variant rotational blur with zero boundary conditions. Both yield the matrix $A\in\R^{256^2 \times 256^2}$. For test image consider the \textit{dotk} image in Figure~\ref{FigureBlurSolutions} (right) polluted by additive Gaussian white noise with 
$\eta=5e-{2}$.

\begin{figure}[ht]
	\centering	
	\includegraphics[width=6cm]{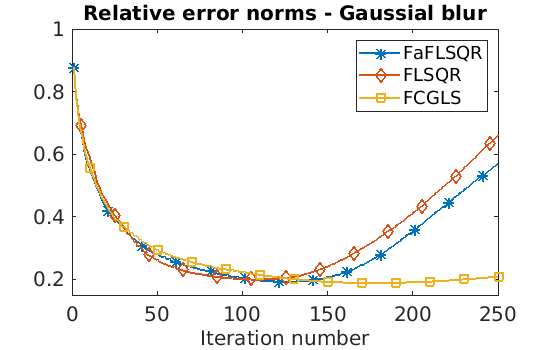}
	\includegraphics[width=6cm]{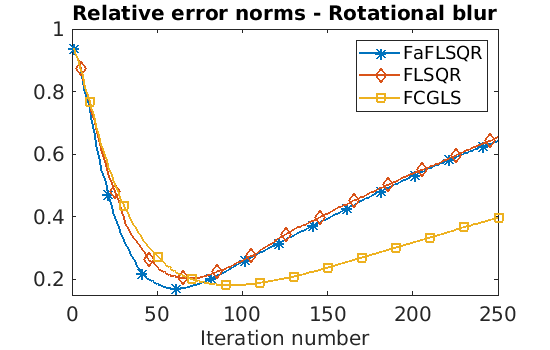}
	\label{FigureBlurError}
	\caption{Experiment 2. Relative error norms with respect to the number of iterations $k$, for FaFLSQR, FLSQR and FCGLS 
	on two setting of the test problem \textit{PRBlur}, $\eta=5e-{2}$.
	}
\end{figure}

Figure~\ref{FigureBlurError} shows the relative error norms. Convergence of FaLSQR and FLSQR is comparable, FaFLSQR attains a slightly lower error norm here. For both test examples, FCGLS exhibits delay of convergence compared to FLSQR and FaFLSQR. The values of the lowest attained error norms are summarized in Table~\ref{TableBlurLowestError}, column 4.

\begin{figure}[ht]
	\centering	
	\includegraphics[width=6cm]{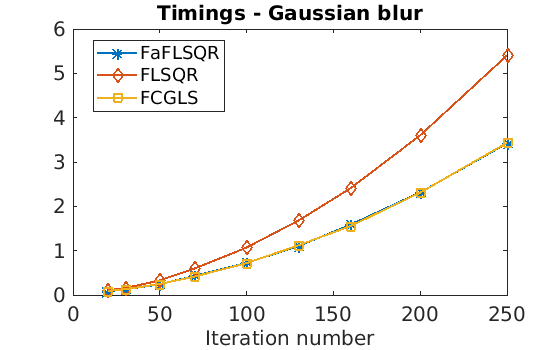}
	\includegraphics[width=6cm]{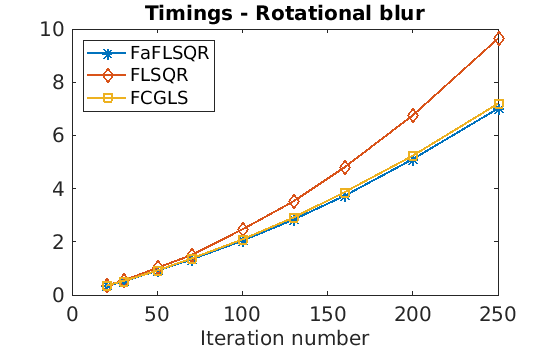}
	\label{FigureBlurTiming}
	\caption{Experiment 2. Computational time (in seconds) with respect to the number of iterations $k$, for FaFLSQR, FLSQR and FCGLS 
	on two setting of the test problem \textit{PRBlur}, $\eta=5e-{2}$.}
\end{figure}  

Figure~\ref{FigureBlurTiming} depicts computational time. FaFLSQR and FCGLS need almost identical time to compute the same number of iterations $k$. FLSQR is slower, the gap increases with $k$ and it is more pronounced for the Gaussian blur with periodic boundary conditions. This is due to the fact that the matrix-vector products with matrices representing point spread functions can be implemented very efficiently (using Kronocker products and Fast fourier transform). Thus the additional long reccurence in FLSQR plays more significant role in the overall computational time. 
Table~\ref{TableBlurLowestError} compares the number of iterations and time required by each method to reach the lowest error norm. Since the convergence is slightly different for each method, Table~\ref{TableBlurSameError} further compares the number of iterations and time required to reach the same error norm for all methods. The speed up value $\tau$ indicates that FaFLSQR requires significantly lower computational time in all considered cases.

Figure~\ref{FigureBlurSolutions} then depicts the best approximate solutions computed by each method. We do not observe significant difference in the quality of the approximations for the Gaussian blur, the lowest relative error norm was very similar for all methods. 
For the Rotational blur, the approximation from FLSQR is slightly more blurred which corresponds to the higher value of relative error norm here.

\begin{table}[htbp]
	\footnotesize
	\caption{Experiment 2. Computational time required by FaLSQR, FLSQR and FCGLS, to find the approximation with the lowest relative error norm. The speed up value $\tau$ of FaFLSQR is in the last column.}\label{TableBlurLowestError}
	\begin{center}
		\begin{tabular}{|c|c|c|c|c|c|c|} \hline
			Blur	&  Method &Iterations & Relative Error & Time (sec) & $\tau$ \\ \hline
			Gaussian &  FaFLSQR & 128 &  0.1941 & 1.26 &  \\
			&  FLSQR & 110 &  0.2009 & 1.47 & 14.2 \% \\
			&  FCGLS & 179 & 0.1901 & 2.31 &  45.5 \% \\ \hline
			Rotational &  FaFLSQR &60 & 0.1717 & 1.11      &  \\
			&  FLSQR &70 & 0.2045 & 1.56 & 28.8 \% \\ 
			&  FCGLS &93 & 0.1824 & 1.93 &  42.5 \% \\ \hline
		\end{tabular}
	\end{center}	
\end{table}

\begin{table}[htbp]
	\footnotesize
	\caption{Experiment 2. Computational time required by FaLSQR, FLSQR and FCGLS, to achieve approximately the same relative error norm. The speed up value $\tau$ of FaFLSQR is in the last column.}\label{TableBlurSameError}
	\begin{center}
		\begin{tabular}{|c|c|c|c|c|c|c|} \hline
			Blur	&  Method &Iterations & Relative Error & Time (sec) & $\tau$ \\ \hline
			Gaussian &  FaFLSQR & 106 & 0.2007 & 0.94 &   \\
			&  FLSQR & 110 &  0.2009 & 1.46 & 35.6 \% \\
			&  FCGLS & 133 & 0.2014 & 1.39 &  32.3 \% \\ \hline
			Rotational &  FaFLSQR &44 & 0.2055	 & 0.79       &  \\
			&  FLSQR &70 & 0.2045 & 1.56 & 49.4 \%  \\ 
			&  FCGLS &125 & 0.2048 & 2.83 & 72.1 \%  \\ \hline
		\end{tabular}
	\end{center}	
\end{table}

\subsection{Experiment 3}\label{Experiment3}
Now we focus closely on FaLSQR and FLSQR. Consider the two dimensional tomography problem \textit{PRtomo} from the IRtools toolbox \cite{IRTools} with $N=256$, $d=\sqrt{2}N$, $angles=0:2:179$, and two values of $p$ (number of rays). The first value $p=round(\sqrt{2}N)$ yields an underdetermined matrix $A\in \R^{32 \, 580 \times 65 \, 536 }$, the second $p=3.5N$ yields an overdetermined matrix $A\in \R^{80 \, 640 \times 65\, 536 }$. For both problems, two levels of additive Gaussian white noise in the right hand side are considered, $\eta_1=1e-{2}$ and  $\eta_2=1e-{3}$. 

\begin{figure}[ht]
	\centering	
	\includegraphics[width=6cm]{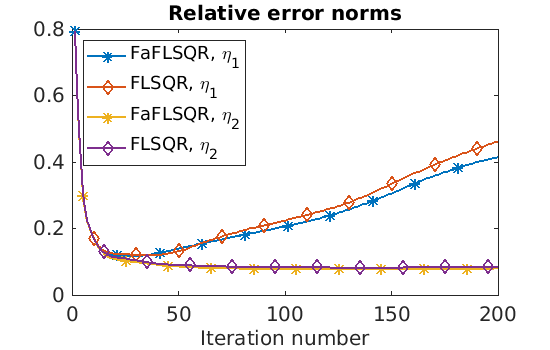}
	\includegraphics[width=6cm]{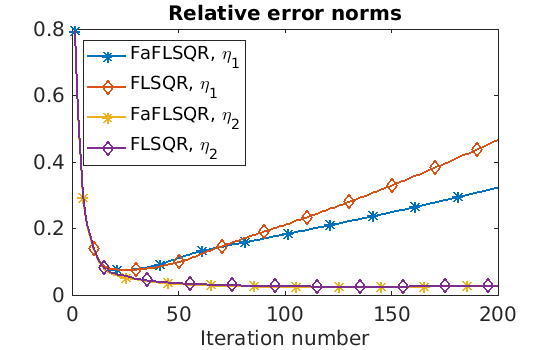}
		\label{FigureTomoError}
	\caption{Experiment 3. Relative error norms with respect to the number of iterations $k$ for FaFLSQR and FLSQR, 
	noise levels $\eta_1=1e-{2}$ and  $\eta_2=1e-{3}$. Test problem
	\textit{PRtomo} with underdetermined (left) and overdetermined (right) matrix, respectively.}
\end{figure}

\begin{figure}[ht]
	\centering	
	\includegraphics[width=6cm]{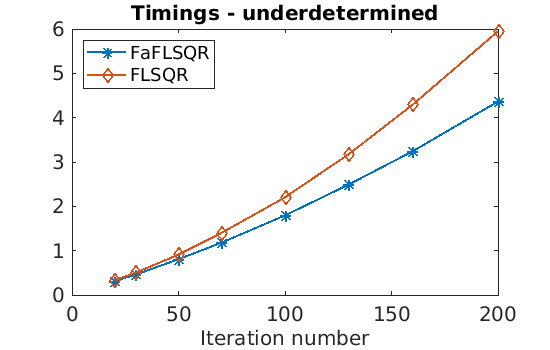}
	\includegraphics[width=6cm]{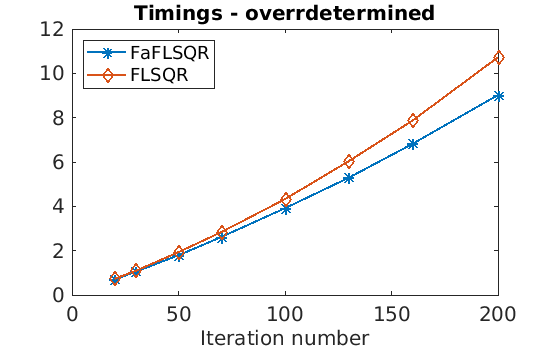}
	\label{FigureTomoTimings}
	\caption{Experiment 3.  Computational time with respect to the number of iterations $k$ for FaFLSQR and FLSQR, 
	noise levels $\eta_1=1e-{2}$ and  $\eta_2=1e-{3}$. Test problem
	\textit{PRtomo} with underdetermined (left) and overdetermined (right) matrix, respectively.}
\end{figure}

Figure~\ref{FigureTomoError} and Figure~\ref{FigureTomoTimings} show the relative error norm  and computational time with respect to number of iterations $k$ for FaLSQR and FLSQR and all four test problems. For the lower level of noise, the error norms behave similarly. For the higher level of noise, the curves are close until the point of semiconvergence, where they start to differ. However, the computational time is higher for FLSQR, see Figure~\ref{FigureTomoTimings}. The difference increases with the increasing $k$ and it is more pronounced for the underdetermined system. This is because the extra full recurrence in FLSQR is computed on the set of longer basis vectors.


\begin{figure}[ht]
	\centering	
	\includegraphics[width=13cm]{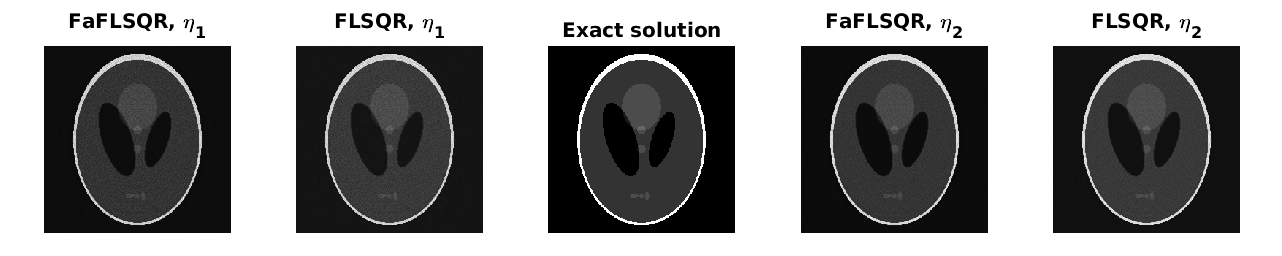}
	\label{FigureTomoSolutions}
	\caption{Experiment 3. Comparison of the exact solution and approximations with the lowest error norms obtained 
	by FaLSQR and FLSQR. Test problem \textit{PRtomo} with underdetermined matrix, noise levels $\eta_1=1e-{2}$ and $\eta_2=1e-{3}$.}
\end{figure}

Table~\ref{TablePRtomo} summarizes the lowest achieved relative error norms and the corresponding computational times of FaFLSQR and FLSQR. 
Despite the number of iterations is very low compared to the size of the problem, we observe FaFLSQR is significantly faster than FLSQR. The lowest relative error norms are similar for all considered settings. Thus, the computed approximate solutions are also close, see Figure~\ref{FigureTomoSolutions} for illustration on the underdetermined system.

\begin{table}[htbp]
	\footnotesize
	\caption{Experiment 3. Computational time required by FaLSQR and FLSQR to find the approximation with the lowest relative error norm. The speed up value $\tau$ of FaFLSQR is in the last column.}\label{TablePRtomo}
	\begin{center}
		\begin{tabular}{|c|c|c|c|c|c|c|} \hline
		System	& Noise Level & Method &Iterations & Error & Time (sec) & $\tau$ \\ \hline
			 	Underdet.& $\eta_1$& FaFLSQR & 29 & 0.1191 & 0.45 &  \\
	 							& & FLSQR & 36 & 0.1204 & 0.62 & 27.4 \% \\
				 		 & $\eta_2$& FaFLSQR & 113 & 0.0794 & 2.11 &  \\
								 & & FLSQR & 135 & 0.0841 & 3.37 & 37.4 \%  \\ \hline
				Overdet. & $\eta_1$& FaFLSQR &21 & 0.0761 & 0.74      & \\
								 & & FLSQR &27 & 0.0777 & 0.98 &   20.4 \% \\
						 & $\eta_2$& FaFLSQR &120 & 0.0254 & 4.87 &  \\  
								 & & FLSQR &134 & 0.0267 & 6.33 &  23.0 \% \\ \hline
		\end{tabular}
	\end{center}	
\end{table}

\begin{figure}[ht]
	\centering	
		\includegraphics[width=6cm]{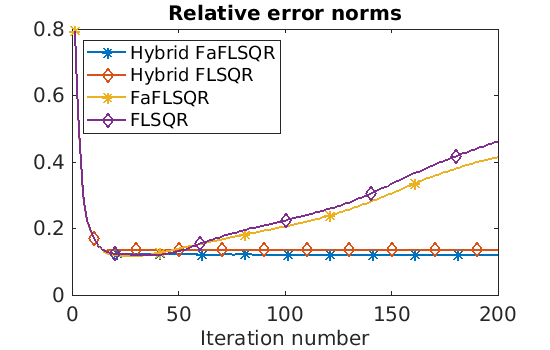}
	\includegraphics[width=6cm]{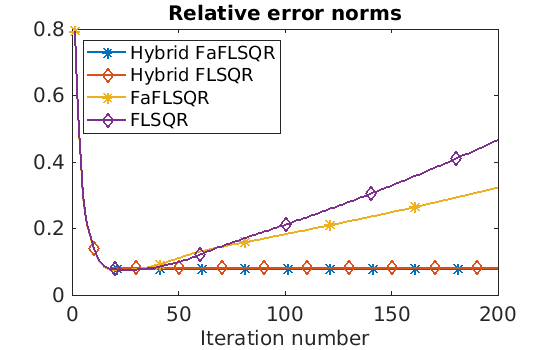}
	\label{FigureTomoHybrid}
	\caption{Experiment 3. Comparison of relative error norms for basic and hybrid variants of FaFLSQR and FLSQR. Test problem \textit{PRtomo} with underdetermined (left) and overdetermined (right) matrix, respectively, noise level $\eta_1=1e-{2}$.}
\end{figure}

Lastly, we compare hybrid variants of FaLSQR and FLSQR on problems with the higher noise level $\eta_1$, where additional regularization might be beneficial. Consider Tikhonov inner regularization of the projected problem in each outer iteration $k$, see \eqref{hybridFLSQR}.
Discrepancy principle is used for the choice of the regularization parameter $\lambda_k$, where we take advantage of the known noise level. Other methods such as GCV, NCP, L-Curve, UPRE, etc., could be applied; see \cite{Hansen} for references. Figure~\ref{FigureTomoHybrid} illustrates that for both methods additional regularization leads to stabilization of the approximation. There is no significant difference in the relative error norms of Hybrid FaFLSQR and Hybrid FLSQR. Thus Hybrid FaFLSQR can be used to solve efficiently problems solvable by Hybrid FLSQR.

\section{Conclusions}
\label{sec:conclusions}
This paper introduced the Fast Flexible LSQR algorithm and its hybrid variant for solving large-scale linear inverse problems. We established theoretical foundations of the method, including the orthogonality of the generated bases vectors, and demonstrated its mathematical equivalence to Flexible CGLS. In contrast to the conventional Flexible LSQR, the proposed approach requires only one long recurrence orthogonalization, leading to a notable reduction in computational cost for sparse or structured systems.
Fast Flexible LSQR achieves greater computational efficiency in floating-point arithmetic than the traditional Flexible LSQR, while maintaining comparable approximation accuracy. Furthermore, the new method exhibits faster convergence than Flexible CGLS due to less significant delay of convergence.

Future research may include alternative preconditioning strategies, e.g., the adaptation of techniques from \cite{nonNegative} to impose nonnegativity constraints in the solution. Another direction is the study of truncated recurrences to further reduce computational overhead. Development of Fast Flexible variant of LSMR algorithm \cite{LSMR} represents a natural extension of this work.

\section{Appendix}

\begin{lemma}\label{LemmaROG}
	The vectors $r_k$ generated by FCGLS given in Algorithm \ref{AlgorithmFCGLS} have the following orthogonality property 
	\[
	r_k \perp q_i,\ i=0,1,\ldots, k-1.
	\]
	\begin{proof}
		We give a proof by induction. Clearly $r_1 \perp q_0$. Assume $r_k\perp q_i$, $i=0,1,\ldots k-1$. From the construction of $r_{k+1}$, 
		obviously $(r_{k+1},q_k)=0$. Moreover, for $i=1,2 \ldots, k-1$, 
		\[
		(r_{k+1},q_i)=(r_k-\gamma_k q_k,q_i)=(r_k,q_i)-\gamma_k(q_k,q_i)=0,
		\]
		using orthogonality of the vectors $q_i$.
	\end{proof}
\end{lemma}
\begin{theorem}
	The vectors $s_k$ generated by FCGLS given in Algorithm \ref{AlgorithmFCGLS} have the following orthogonality property
	\[
	s_k \perp_{M^{-1}_{i+1}} s_i,\ i=0,1,\ldots, k-1.
	\]
	\begin{proof}	
		Let $i, k \in \N$, $0 \leq i< k$. Then
		\begin{align}\label{derivation_sksi}
			\begin{split}
				(s_k,M_{i+1}^{-1}s_{i})
				&=\left(A^Tr_k,M_{i+1}^{-1}s_{i}\right)\\
				&=\left(A^T(r_{k-1}-\gamma_{k-1}q_{k-1}),M_{i+1}^{-1}s_{i}\right)\\
				&=(A^Tr_{k-1},M_{i+1}^{-1}s_i)-\gamma_{k-1}(A^Tq_{k-1},\hat{s}_i)\\
				&=(s_{k-1},M_{i+1}^{-1}s_i)-\gamma_{k-1}(q_{k-1},A\hat{s}_i)\\
				&=(s_{k-1},M_{i+1}^{-1}s_i)-\gamma_{k-1}\left(q_{k-1},q_i-\sum_{j=0}^{i-1}\theta^{i-1}_{j}q_j \right).
			\end{split}
		\end{align}	
		For $i=(k-1)$, we have 
		\begin{align*}
			(s_{k-1},M_{k}^{-1}s_{k-1})&=(A^Tr_{k-1},\hat{s}_{k-1})=(r_{k-1},A\hat{s}_{k-1})\\
			&=\left(r_{k-1},q_{k-1}-\sum_{j=0}^{k-2}\theta^{k-2}_{j}q_j\right)= (r_{k-1},q_{k-1}) \\
			\gamma_{k-1}\left(q_{k-1},q_{k-1}-\sum_{j=0}^{k-2}\theta^{k-2}_{j}q_j \right)
			&=\left( (r_{k-1},q_{k-1})/\norm{q_{k-1}}^2 \right)(q_{k-1},q_{k-1})\\
			&= (r_{k-1},q_{k-1}),
		\end{align*}
		using the orthogonality of the vectors $q_i$ and Lemma \ref{LemmaROG}.					
		For $i<(k-1)$, we have 
		\begin{equation*}
			\gamma_{k-1}\left(q_{k-1},q_i-\sum_{j=0}^{i-1}\theta^{i-1}_{j}q_j \right)=0.
		\end{equation*}
		Substitution into \eqref{derivation_sksi} gives
		\begin{align}
			(s_k,M_{k}^{-1}s_{k-1})&= 0 \\
			(s_k,M^{-1}_{i+1}s_{i})&=(s_{k-1},M^{-1}_{i+1}s_{i}),\ i=0, \ldots, k-2. \label{expression_sksi}
		\end{align}			
		Now proceed by induction. For $k=1$, clearly $s_1\perp_{M_{1}^{-1}}s_0$. Assume that 
		$s_{k-1} \perp_{M^{-1}_i} s_i,\ i=0,1,\ldots, k-2$. Then \eqref{expression_sksi} immediately gives 
		\[
		s_k\perp_{M^{-1}_{i+1}}s_i, \ i=0, 1, \ldots k-1.
		\]
	\end{proof}	
\end{theorem}

\section*{Acknowledgments}

\bibliographystyle{siamplain}
\bibliography{references}
\end{document}